\documentclass[12pt]{article}
\begin{document}

\begin{center}
{\Large Partial Inertial Manifolds for infinite-dimensional dynamical
systems: Example for P.D.E.s with a state-dependent delay}

\bigskip

{\sc Alexander V. Rezounenko}

\smallskip

 Department of Mechanics and Mathematics, Kharkov University,

 4, Svobody Sqr., Kharkov, 61077, Ukraine

 E-mail: rezounenko@univer.kharkov.ua

\end{center}

\begin{quote}
{\bf Abstract.} We propose a new notion of Partial Inertial
Manifold to study the long-time asymptotic behavior of dissipative
differential equations. As shown on an example, such manifolds may
exist in the cases when the classical Inertial manifold does not
exist (or not known to exist).
\end{quote}




\bigskip
{\bf 1. Introduction} 

\medskip

Study of the long-time asymptotic behavior of solutions occupies
an important place in the qualitative theory of differential
equations. Considering partial and/or functional differential
equations one naturally obtains infinite-dimensional dynamical
systems. To investigate their asymptotic behavior many powerful
methods and approaches have been developed, such as global, weak
and exponential attractors
\cite{Babin-Vishik,efnt,Temam_book,Chueshov_book}, inertial
manifolds \cite{fstem,fstiti,Temam_book,chowlu,con-f-n-tem},
approximate inertial
manifolds \cite{fmt,fstem,
Temam_book}, determining functional \cite{Chueshov_book} etc.

During these investigations many deep results were obtained so far
and the subject continuously attracts attention of many
researchers. Each of the mentioned objects (attractors, manifolds,
functionals) indicates important features of the dynamical systems
under considerations, but naturally has special conditions to
exist. If we are able to establish simultaneously the existence of
several of the mentioned objects for a system, then we get more
important information on its asymptotic properties. In this note
we introduce a new notion - Partial Inertial Manifold and hope it
will be useful for the study.

\medskip
{\bf 2. Partial Inertial Manifolds}

\medskip

Consider a dynamical system $(S(t), {\cal H}),$ where $S(t) : {\cal H} \to
{\cal H}$ denotes the evolution operator and ${\cal H}$ is the phase space
(see e.g.
\cite{Babin-Vishik,Hale_book,Temam_book,Walther_book,Chueshov_book} for
more details). For example, one may consider a general dissipative
differential equation in the space ${\cal H}$
\begin{equation}\label{sdd6-1}
\dot u+Au =B(u), \quad u\in {\cal H},
\end{equation}
where $A$ is the (leading in some sense) linear part, and $B$ is the
nonlinearity. Under the  natural assumptions this equation generates an
evolution operator as a shift along the trajectories of (\ref{sdd6-1})
i.e. $S(t) u^0\equiv u(t;u^0),$ where $u(t;u^0)$ denotes the solution of
(\ref{sdd6-1}) with the initial data $u(0)=u^0.$

Such objects as global attractors and inertial manifolds play an
important role in the study of long-time asymptotic behavior of
dissipative dynamical systems. We recall
\cite{fstem,fstiti,Temam_book,chowlu,con-f-n-tem}

\medskip

{\bf Definition~1.} {\it A set ${\cal M}\subset {\cal H}$ is called an
{\rm Inertial manifold} if there exist a projector $P=P^2: {\cal H}\to
{\cal H}$ and a Lipschitz mapping $\Phi: P{\cal H}\to (1-P){\cal H}$ such
that
\begin{itemize}
 \item $\dim P <\infty;$
 \item ${\cal M}=\left\{\, u : u=p+\Phi(p),\quad p\in P{\cal H}\, \right\} \subset {\cal H};$
 \item $S(t) {\cal M} \subset {\cal M}$ for all $t\ge 0;$
 \item for any $u\in {\cal H}$ one has $\hbox{ dist}_{\cal H} \{ S(t) u, {\cal
 M}\}\le K(||u||_{\cal H})\cdot \exp \{ -\alpha t\} $ for some $\alpha>0.$
\end{itemize} }

 The existing theory says that a dynamical system usually has an
Inertial manifold provided special {\it spectral gap conditions} are
satisfied (see e.g. \cite{Babin-Vishik,Temam_book,Chueshov_book} for more
details). These conditions are usually formulated as a condition for the
distance between two nearest eigenvalues $|\lambda_{N+1}-\lambda_N|$ of
the leading linear part of the differential equation to be big enough in
comparison with the Lipschitz constant of the nonlinear part of the
differential equation and (possibly) lower degrees $\lambda^\alpha_{N+1},
\lambda^\alpha_{N}, \alpha\in [0,1)$ of the eigenvalues. In this
direction, to get an inertial manifold, one first computes the Lipschitz
constant $L$ of the nonlinear part $B$ and than looks for an integer $N$
such that
$|\lambda_{N+1}-\lambda_N|\ge C (L,\lambda^\alpha_{N+1}, \lambda^\alpha_{N})$ 
(to be more precise, one needs to consider a concrete equation).
Unfortunately, the spectral gap conditions are very restrictive and do not
hold for many important problems. To investigate the cases when inertial
manifold does not exist (or not known to exist) 
another approaches have been proposed such as approximate inertial
manifolds, exponential attractors etc (see e.g.
\cite{Temam_book,Chueshov_book}).

In this note we propose a new approach. The main idea is to look for a
subset $D$ of the phase space ${\cal H}$ such that the restriction of the
nonlinear term of the differential equation on the set $\bigcup_{t\ge
0}S(t) D\subset {\cal H}$ has a small enough Lipschitz constant. If we are
able to extend  the restriction of the nonlinear term from $\bigcup_{t\ge
0}S(t) D$ to ${\cal H}$ without increasing the Lipschitz constant, then we
get an auxiliary nonlinear term $B_\ell.$ If the spectral gap conditions
are satisfied with this (smaller) Lipschitz constant, then equation
(\ref{sdd6-1}) with the nonlinearity $B_\ell$ does have an inertial
manifold. This manifold is finite-dimensional and attracts all the
trajectories of the initial equation (\ref{sdd6-1}) which start in
$\bigcup_{t\ge 0}S(t) D.$ We call this manifold {\it partial inertial
manifold} for (\ref{sdd6-1}). The name reflects the fact that the manifold
attracts only {\it part} of the phase space, but not the whole ${\cal H}.$
Considerations become simpler if the set $D$ is positively invariant i.e.
$S(t) D\subset D$ for all $t\ge 0,$ then  $\bigcup_{t\ge 0}S(t) D=D.$

We summarize the above ideas in the following
\medskip

{\bf Definition~2.} {\it A set ${\cal M}\subset {\cal H}$ is called a {\rm
Partial Inertial Manifold} if there exist a projector $P=P^2: {\cal H}\to
{\cal H}$, a Lipschitz mapping $\Phi: P{\cal H}\to (1-P){\cal H}$ and a
set $D\subset {\cal H}$ such that
\begin{itemize}
 \item $\dim P <\infty;$
 \item ${\cal M}=\left\{\, u : u=p+\Phi(p),\quad p\in P{\cal H}\, \right\} \subset {\cal H};$
 \item for any $u\in D\subset {\cal H}$ one has $\hbox{ dist}_{\cal H} \{ S(t) u, {\cal
 M}\}\le K(||u||_{\cal H})\cdot \exp \{ -\alpha t\} $ for some $\alpha>0.$
\end{itemize} }
\medskip

{\bf Remark.} {\it It is easy to see that Definition~2 gives the
possibility to exist more than one Partial Inertial Manifolds for
the same equation if we have several sets $D_i$ with the described
properties. On the other hand, the classical Inertial Manifold is
a Partial Inertial Manifold if we set $D={\cal H}$.}

\medskip


In the next section we present a concrete example of a system of
partial differential equations with state-dependent distributed
delay for which a partial inertial manifold exists while inertial
manifold does not. The construction of the example is based on our
recent studying of P.D.E.s with state-dependent delay
\cite{Rezounenko-Wu-2006,Rezounenko-JMAA-2007}. For more details
on state-dependent (ordinary) equations see e.g.
\cite{MalletParet,Walther_JDE-2003}.

\bigskip {\bf 3. Example of the existence of a P.I.M.: state-dependent delay equations}

\medskip

Consider the following 
partial differential equation with state-dependent distributed delay
\begin{equation}\label{sdd6-g}
 \frac{\partial }{\partial t}u(t,x)+Au(t,x)  
= \int^0_{-r} b(u(t+\theta, x)) \xi(\theta, u_t) d\theta \equiv \big(
B_1[\xi](u_t) \big)(x),\quad x\in \Omega,
\end{equation}
 where $A$ is a densely-defined self-adjoint positive linear operator
 with domain $D(A)\subset L^2(\Omega )$ and with compact
  resolvent, so $A: D(A)\to L^2(\Omega )$ generates an analytic semigroup,
  $\Omega $ is a smooth bounded domain in $R^{n_0}$,
  $b:R\to R$ is a locally Lipschitz bounded map
  ($|b(w)|\le M_b$ with $M_b\ge 0),$ 
  The function $\xi (\cdot,\cdot): [-r,0]\times C
 \to R$ represents the state-dependent  distributed delay.
   We denote for short $C\equiv
C([-r,0]; L^2(\Omega)).$ As usually for delay equations, we denote
by $u_t$ the function of $\theta\in [-r,0]$ by the formula
$u_t\equiv u_t(\theta)\equiv u(t+\theta).$ For more details on
delay equations we refer to the classical monographs
\cite{Hale_book,Walther_book,Wu_book,Krisztin-Walther-Wu}.

We consider equation (\ref{sdd6-g}) with the following initial conditions
\begin{equation}\label{sdd6-ic}
    u|_{[-r,0]}=\varphi \in C\equiv C([-r,0]; L^2(\Omega)).
\end{equation}

The methods used in our work can be applied to another types of
nonlinear and delay PDEs. We choose a particular form of nonlinear
delay term $B_1$ for simplicity and to illustrate our approach on
the diffusive Nicholson's blowflies equation (see below for more
details).

Assume the following:
\begin{equation}\label{sdd6-b1}
{\bf A1)}\quad |b(s)|\le M_b \hbox{ and } |b(s^1)-b(s^2)|\le L_b|s^1-s^2|,
\,\,\mbox{for all} \,\, s,s^1,s^2\in R.
\end{equation}
\begin{equation}\label{sdd6-xi1}
{\bf A2)}\quad \int^0_{-r}|\xi(\theta, \psi^1)-\xi (\theta, \psi^2) | d\theta\\
\le L^{1,1}_{\xi,M}
  \cdot  ||\psi^1-\psi^2||_{L^1(-r,0;L^1(\Omega))},
\end{equation}
\begin{equation}\label{sdd6-xi2}
{\bf A3)}\qquad  \qquad  \qquad  ess \sup_{\theta\in (-r,0)}
|\xi(\cdot,\psi)|\le M_\xi \,\,\mbox{for all} \,\,  \psi\in C.
\end{equation}
\medskip

We notice that assumptions (\ref{sdd6-b1})-(\ref{sdd6-xi2}) are more
restrictive than the ones of \cite[theorems~1,2]{Rezounenko-JMAA-2007}, so
we can apply theorems~1,2 from \cite{Rezounenko-JMAA-2007} to get the
existence and uniqueness of solutions for (\ref{sdd6-g}), (\ref{sdd6-ic})
with $\varphi\in C.$ In this note we are interested in continuous
solutions i.e. functions $u\in C([0,T]; L^2(\Omega))$ for any $T>0.$

In the same manner, using
\cite[theorems~1,2]{Rezounenko-JMAA-2007}, we define an {\it
evolution operator} $S_t : C\to C$ by the formula $S_t
=u_t(\varphi),$ where  $u(\varphi)$ denotes the unique
(continuous) solution of (\ref{sdd6-g}),(\ref{sdd6-ic}) with the
initial condition $u_0(\varphi)=\varphi.$ Sometimes, we will write
$S_t[\xi]$ to indicate the kernel function $\xi$ in the
nonlinearity $B_1[\xi]$ (see (\ref{sdd6-g})).

\medskip

Notice that due to the the inclusion $C \subset L^1(-r,0; L^1(\Omega)) ,$
we get for any $v^1,v^2\in C:$
$$
||v^1-v^2||_{L^1(-r,0; L^1(\Omega))} = \int^0_{-r} \left\{ \int_\Omega
\left| v^1(\theta,x)-v^2(\theta,x)\right| dx \right\} d\theta
$$
\begin{equation}\label{sdd6-12}
\le \sqrt{|\Omega|} \int^0_{-r}
||v^1(\theta,\cdot)-v^2(\theta,\cdot)||_{L^2(\Omega))}\le r\cdot
\sqrt{|\Omega|}\cdot ||v^1-v^2||_C.
\end{equation}
Hence (\ref{sdd6-xi1}) implies
\begin{equation}\label{sdd6-xi}
\quad \int^0_{-r}|\xi(\theta, \psi^1)-\xi (\theta, \psi^2) | d\theta\\
\le L^{1,1}_{\xi,M}\cdot r\cdot \sqrt{|\Omega|}
  \cdot  ||v^1-v^2||_C.
\end{equation}

Let us check that the mapping  $B_1\equiv B_1[\xi] : C\to L^2(\Omega)$
satisfies the Lipschitz property (c.f. (2.3) in
\cite{Boutet-Chueshov-Rezounenko_NA-1998}). Using (\ref{sdd6-xi}), one has
$$ ||B_1(v^1_0)-B_1(v^2_0)||^2 = \int_\Omega
\left| \int^0_{-r} \left\{ b(v^1(\theta, x)) \xi(\theta,
v^1_0)-b(v^2(\theta, x)) \xi(\theta, v^2_0)\right\} d\theta \right|^2\, dx
$$ $$\le 2L^2_b M^2_\xi \int_\Omega
\left(\int^0_{-r} |v^1(\theta, x)-v^2(\theta, x)| d\theta\right)^2 dx +
2M^2_b \left( L^{1,1}_{\xi,M}\right)^2 r^2 |\Omega|\cdot ||v^1-v^2||^2_C
$$
$$\le 2L^2_b M^2_\xi r \int_\Omega
\left\{\int^0_{-r} |v^1(\theta, x)-v^2(\theta, x)|^2 d\theta\right\} dx +
2M^2_b \left( L^{1,1}_{\xi,M}\right)^2 r^2  |\Omega|\cdot ||v^1-v^2||^2_C
$$
$$\le 2L^2_b M^2_\xi r\int^0_{-r} ||v^1(\theta, \cdot)-v^2(\theta, \cdot)||^2_{L^2(\Omega)}\,
d\theta + 2M^2_b \left( L^{1,1}_{\xi,M}\right)^2 r^2  |\Omega|\cdot
||v^1-v^2||^2_C
$$
$$\le 2\left( L^2_b M^2_\xi r^2  +
M^2_b \left( L^{1,1}_{\xi,M}\right)^2 r^2  |\Omega|\right)\cdot
||v^1-v^2||^2_C.
$$
So, we get (c.f. (2.3) in \cite{Boutet-Chueshov-Rezounenko_NA-1998})
\begin{equation}\label{sdd6-2}
||B_1(v^1)-B_1(v^2)|| \le M_1\cdot ||v^1-v^2||_C \quad \hbox{ with }\quad
M_1\equiv r\cdot \sqrt{2\left( L^2_b M^2_\xi   + M^2_b \left(
L^{1,1}_{\xi,M}\right)^2  |\Omega|\right)}.
\end{equation}
To get estimate (\ref{sdd6-2}), we used the following calculations
$$\left|
\int^0_{-r} \left\{ b(v^1(\theta, x)) \xi(\theta, v^1_0)-b(v^2(\theta, x))
\xi(\theta, v^2_0)\right\} d\theta \right| $$ 
$$\le L_b  \int^0_{-r}
|v^1(\theta, x)-v^2(\theta, x)|\cdot |\xi(\theta, v^1_0)| d\theta + M_b
\int^0_{-r} | \xi(\theta, v^1_0)-\xi(\theta, v^2_0)| d\theta $$ 
$$\le L_b
M_\xi \int^0_{-r} |v^1(\theta, x)-v^2(\theta, x)| d\theta + M_b
L^{1,1}_{\xi,M} \cdot r\cdot \sqrt{|\Omega|} \cdot ||v^1_0- v^2_0)||_C$$
and the inclusion $C\subset L^2((-r,0)\times \Omega) ,$ which implies
$\int^0_{-r} |v^1(\theta, x)-v^2(\theta, x)| d\theta \le
\sqrt{r}\cdot\left( \int^0_{-r} |v^1(\theta, x)-v^2(\theta, x)|^2 d\theta
\right)^{1\over 2}$ and, as a result, $\left(\int^0_{-r} |v^1(\theta,
x)-v^2(\theta, x)| d\theta\right)^2 \le r\cdot \int^0_{-r} |v^1(\theta,
x)-v^2(\theta, x)|^2 d\theta.$

\medskip

Now we recall a sufficient conditions for the existence of an inertial
manifold in the case of delay semilinear parabolic equations
\cite{Boutet-Chueshov-Rezounenko_NA-1998}.

Since $A: D(A)\subset L^2(\Omega )\to L^2(\Omega )$ is a densely-defined
self-adjoint positive linear operator, then there exists an orthonormal
basis $\{ e_k\}$ of $L^2(\Omega )$ such that
$$ Ae_k=\lambda_k e_k, \quad \hbox{ with }\quad 0<\lambda_1\le
\lambda_2\le \ldots, \quad \lim\limits_{k\to\infty} \lambda_k =\infty.$$

As in \cite{Boutet-Chueshov-Rezounenko_NA-1998}, we fix an integer $N$ and
denote $P=P_N$ the orthogonal projector onto the space spanned by the
first $N$ eigenvectors of $A.$ We also define the $N$-dimensional
projector $\hat P= \hat P_N$ in $C$ by
$$\hat P \phi=\left(\hat P \phi\right)(\theta)=\sum^N_{k=1} e^{-\lambda_k \theta}
\langle \phi(0), e_k\rangle_{L^2(\Omega)}\cdot e_k \equiv
e^{-A\theta}\phi(0),\quad \phi\in C, \quad \theta\in [-r,0].
$$

From the above considerations we see that one can apply theorem~3.1 from
\cite{Boutet-Chueshov-Rezounenko_NA-1998} to the system~(\ref{sdd6-g})
under the following assumptions (see
\cite{Boutet-Chueshov-Rezounenko_NA-1998}):

{\bf A4)} For some $N$ and $\mu>0$ the following {\it spectral gap
condition} is satisfied $\lambda_{N+1}-\lambda_N\ge 2\mu $ (see (2.8) in
\cite{Boutet-Chueshov-Rezounenko_NA-1998});

{\bf A5)} Constants $\mu, N$ and delay $r>0$ satisfy: $\mu > 4M_1$ and
$\delta\equiv {2\over \mu} M_1\cdot e^{(\lambda_N+\mu)r}\le {1\over 2}$
(see (3.1) in \cite{Boutet-Chueshov-Rezounenko_NA-1998}).

These two assumptions give (theorem~3.1 from
\cite{Boutet-Chueshov-Rezounenko_NA-1998}) the existence of the
$N$-dimensional asymptotically complete manifold (inertial
manifold) $${\cal M}=\{ \hat p(\theta)+ \Phi (\hat p(0),\theta) :
\hat p(\theta)\in \hat P C\} \subset C$$ which is invariant for
solutions of (\ref{sdd6-g}), (\ref{sdd6-ic}). Here $\Phi$ is a
Liprschitz map $\Phi : P L^2(\Omega) \to (1-\hat P)C.$

If we choose the biggest possible value of constant $\mu={1\over
2}(\lambda_{N+1}-\lambda_N),$ then we get an estimate for the upper bound
of the Lipschitz constant $M_1$:
\begin{equation}\label{sdd6-3}
M_1 \le {\lambda_{N+1}-\lambda_N\over 8}\cdot
\exp\left\{-{(\lambda_{N+1}+\lambda_N)\over 2}\cdot r\right\}.
\end{equation}
Our goal is to illustrate that in the case when (\ref{sdd6-3}) does not
hold, it is possible that partial inertial manifolds do exist.

\bigskip  %

{\bf 3.1. Construction of the kernel function $\xi$.}

\medskip


Let us choose
\begin{equation}\label{sdd6-4}
\xi^{+}(\theta)\ge 0 \hbox{ a. e. in } \theta\in (-r,0)\quad \hbox{ and }
\quad \xi^{-}(\theta)\le 0 \hbox{ a. e. in } \theta\in (-r,0)
\end{equation}
such that \begin{equation}\label{sdd6-5} ess \sup_{\theta\in
(-r,0)}|\xi^{\pm}(\theta)| \le {1\over 2}\, M_\xi.
\end{equation}
For any $v\in C$ we write
\begin{equation}\label{sdd6-6}
v(\theta,x) = v^{+}(\theta,x) + v^{-}(\theta,x)
\end{equation}
where (a.e. in $x\in \Omega$)
\begin{equation}\label{sdd6-7}
v^{+}(\theta,x)\equiv \sup \{ v(\theta,x), 0\} \ge 0, \qquad
v^{-}(\theta,x)\equiv \inf \{ v(\theta,x), 0\} \le 0.
\end{equation}
We will use the following property
\begin{equation}\label{sdd6-8}
||v||_{L^1(-r,0;L^1(\Omega))} = ||v^{+}||_{L^1(-r,0;L^1(\Omega))}
+||v^{-}||_{L^1(-r,0;L^1(\Omega))}.
\end{equation}
Now we are ready to define for any $v\in C$
\begin{equation}\label{sdd6-9}
\xi(\theta, v)=  \xi^{+}(\theta)\cdot
\min\left\{||v^{+}||_{L^1(-r,0;L^1(\Omega))},1\right\} +
  \xi^{-}(\theta)\cdot \min\left\{||v^{-}||_{L^1(-r,0;L^1(\Omega))},1\right\}.
\end{equation}
Using the property (for any norm $||\cdot ||$)
\begin{equation}\label{sdd6-10}
\min\left\{||\psi^1||,1\right\}-\min\left\{||\psi^2||,1\right\}\le
||\psi^1|| - ||\psi^2|| \le ||\psi^1-\psi^2||,
\end{equation}
one can check that $\xi,$ defined in (\ref{sdd6-9}), satisfies 
(\ref{sdd6-xi1}) with
\begin{equation}\label{sdd6-11}
L^{1,1}_{\xi,M}\equiv \max \left\{ \int^0_{-r} |\xi^{+}(\theta)|d\theta,
\int^0_{-r} |\xi^{-}(\theta)|d\theta\right\}.
\end{equation}
More precisely (we will write $||\cdot||_{L^{1,1}}\equiv
 ||\cdot||_{L^1(-r,0;L^1(\Omega))}$ for short):
$$\int^0_{-r}
|\xi(\theta,v^1)-\xi(\theta,v^2)|d\theta = \int^0_{-r} |
\xi^{+}(\theta)\cdot
\left[\min\left\{||v^{1+}||_{L^{1,1}},1\right\}-\min\left\{||v^{2+}||_{L^{1,1}},1\right\}\right]
$$
$$+
  \xi^{-}(\theta)\cdot \left[\min\{||v^{1-}||_{L^{1,1}},1\}-\min\{||v^{2-}||_{L^{1,1}},1\}
  \right]|d\theta
$$
$$ \le \int^0_{-r}\left\{ \left|\, \xi^{+}(\theta)\right| \cdot ||v^{1+}-v^{2+}||_{L^{1,1}}  +
\left|\, \xi^{-}(\theta)\right| \cdot ||v^{1-}-v^{2-}||_{L^{1,1}} \right\}
d\theta
$$
$$\le \max \left\{ \int^0_{-r} |\xi^{+}(\theta)|d\theta,
\int^0_{-r} |\xi^{-}(\theta)|d\theta\right\} \cdot \left(
||v^{1+}-v^{2+}||_{L^{1,1}} + ||v^{1-}-v^{2-}||_{L^{1,1}}  \right)
$$
$$\le L^{1,1}_{\xi,M} \cdot ||v^{1}-v^{2}||_{L^{1,1}},
$$
where  $L^{1,1}_{\xi,M}$ is defined by (\ref{sdd6-11}). Here we also use
(\ref{sdd6-8}).

Definition (\ref{sdd6-9}) and assumption (\ref{sdd6-5}) give
(\ref{sdd6-xi2}). Hence we conclude that function $\xi$, defined by
(\ref{sdd6-9}), satisfies assumptions (\ref{sdd6-xi1}), (\ref{sdd6-xi2}).

\bigskip  %

{\bf 3.2. Properties of the delay term $B_1[\xi]$.}

\medskip


Let us define $D_{+}\equiv \left\{ v\in C : \forall \theta\in [-r,0]
\Rightarrow v(\theta,x) \ge 0 \quad\hbox{ a. e. in } x\in\Omega
\right\}\subset C$ and $D_{-}\equiv \left\{ v\in C : \forall \theta\in
[-r,0] \Rightarrow v(\theta,x) \le 0 \quad\hbox{ a. e. in } x\in\Omega
\right\}\subset C.$

In addition to (\ref{sdd6-b1}), we assume that function $b$ satisfies
\begin{equation}\label{sdd6-b2}
{\bf A6)}\quad  b(s)=b(-s)\ge 0, \quad s\in R.
\end{equation}
So definitions (\ref{sdd6-9}), (\ref{sdd6-4}) and assumption
(\ref{sdd6-b2}) give
$$\forall v\in D_{+} \Rightarrow B_1(v)\ge 0 \quad\hbox{ a. e. in }
x\in\Omega, \quad\hbox{ and }\quad \forall v\in D_{-} \Rightarrow
B_1(v)\le 0 \quad\hbox{ a. e. in } x\in\Omega.
$$
The last property implies (see \cite{henry}) that cones $D_{+},D_{-}$ are
positively invariant i.e.
\begin{equation}\label{sdd6-13}
S_t[\xi] D_{+} \subset D_{+} \quad\hbox{ and }\quad S_t[\xi] D_{-} \subset
D_{-}.
\end{equation}
Here $S_t[\xi]: C\to C$ denotes the evolution operator constructed by the
solutions of (\ref{sdd6-g}), (\ref{sdd6-ic}) with the kernel function
$\xi$ in (\ref{sdd6-g}), defined by (\ref{sdd6-9}).

Now we consider two auxiliary functions (see (\ref{sdd6-9}),
(\ref{sdd6-6}), (\ref{sdd6-7}))
\begin{equation}\label{sdd6-14}
\xi^p(\theta, v)\equiv  \xi^{+}(\theta)\cdot
\min\left\{||v^{+}||_{L^1(-r,0;L^1(\Omega))},1\right\},
\end{equation}
\begin{equation}\label{sdd6-15}
 \xi^n(\theta, v)\equiv
  \xi^{-}(\theta)\cdot \min\left\{||v^{-}||_{L^1(-r,0;L^1(\Omega))},1\right\}.
\end{equation}

Since $\forall v\in D_{+} \Rightarrow \xi^p(\theta, v)=\xi(\theta, v),$
then property (\ref{sdd6-13}) gives
\begin{equation}\label{sdd6-16}
\forall v\in D_{+} \Rightarrow S_t[\xi^p] v=S_t[\xi] v.
\end{equation}
 In the same way, $\forall v\in D_{-} \Rightarrow S_t[\xi^n] v=S_t[\xi] v.$

The above considerations clearly show that $B_1[\xi^p]$ satisfies
(\ref{sdd6-xi1}) with the Lipschitz constant $M_1=M_1[\xi^p]$ defined by
(\ref{sdd6-2}) where the constant $L^{1,1}_{\xi^p,M}=\int^0_{-r}
|\xi^{+}(\theta)|d\theta$ instead of $L^{1,1}_{\xi,M}=\max \left\{
\int^0_{-r} |\xi^{+}(\theta)|d\theta, \int^0_{-r}
|\xi^{-}(\theta)|d\theta\right\}$ (see (\ref{sdd6-11})). In the same
manner, we get the Lipschitz constant for $B_1[\xi^n]$ by (\ref{sdd6-2})
with $L^{1,1}_{\xi^n,M}=\int^0_{-r} |\xi^{-}(\theta)|d\theta.$

Due to the explicit dependence of the Lipschitz constants $M_1=M_1[\xi^p]$
and $M_1=M_1[\xi^n]$ on the values $\int^0_{-r} |\xi^{+}(\theta)|d\theta,$
$\int^0_{-r} |\xi^{-}(\theta)|d\theta$ (see (\ref{sdd6-2})), we may choose
small enough value of $\int^0_{-r} |\xi^{+}(\theta)|d\theta$ and big
enough value of $\int^0_{-r} |\xi^{-}(\theta)|d\theta$ such that the
constant $M_1[\xi^p]$ satisfies (\ref{sdd6-3}) while $M_1[\xi^n]$ does
not. Of course, we also need the value $rL_bM_\xi$ to be small enough (see
(\ref{sdd6-2})). In this case, by (\ref{sdd6-11}), the constant $M_1[\xi]$
does not satisfy (\ref{sdd6-3}).

\medskip

{\bf Remark.} {\it More precisely, Let us first choose and fix $r$
small enough to satisfy (see (\ref{sdd6-3}))
\begin{equation}\label{sdd6-17}
r\le  {\lambda_{N+1}-\lambda_N \over 16 L_b M_\xi}
\cdot\exp\left\{ - {\lambda_{N+1}+\lambda_N\over 2}\cdot
r\right\}.
\end{equation}
Then, for the fixed value of $r$, choose $\xi^{+}(\cdot)$ such
that
\begin{equation}\label{sdd6-18}
\int^0_{-r}|\xi^{+}(\theta)|\, d\theta\le {\lambda_{N+1}-\lambda_N
\over 16\, r M_b  \sqrt{|\Omega|}} \cdot\exp\left\{ -
{\lambda_{N+1}+\lambda_N\over 2}\cdot r\right\}.
\end{equation}
Assumptions (\ref{sdd6-17}), (\ref{sdd6-18}) imply that
$M_1[\xi^p]$ satisfies (\ref{sdd6-3}). Now we choose
$\xi^{-}(\cdot)$ such that
\begin{equation}\label{sdd6-19}
\int^0_{-r}|\xi^{-}(\theta)|\, d\theta > {\lambda_{N+1}-\lambda_N
\over 8\, r M_b  \sqrt{|\Omega|}} \cdot\exp\left\{ -
{\lambda_{N+1}+\lambda_N\over 2}\cdot r\right\}.
\end{equation}
Assumptions (\ref{sdd6-17}), (\ref{sdd6-19}) imply that $M_1[\xi]$
and $M_1[\xi^n]$ do not satisfy (\ref{sdd6-3}).
}%

\smallskip

These considerations clearly show that the system (\ref{sdd6-g}),
(\ref{sdd6-ic}) with the right hand side $B_1[\xi]$ ($\xi$ defined by
(\ref{sdd6-9})) does not possess an inertial manifold, while the system
(\ref{sdd6-g}), (\ref{sdd6-ic}) with the right hand side $B_1[\xi^p]$
($\xi^p$ defined by (\ref{sdd6-14})) does possess (due to
\cite[theorem~3.1]{Boutet-Chueshov-Rezounenko_NA-1998}). Since the
evolution operators $S_t[\xi^p]$ and $S_t[\xi] $ coincide on $D_{+}$ (see
(\ref{sdd6-16})), we may conclude that the system (\ref{sdd6-g}),
(\ref{sdd6-ic}) with the right hand side $B_1[\xi]$ ($\xi$ defined by
(\ref{sdd6-9})) possesses a finite-dimensional manifold (inertial manifold
for the system with $B_1[\xi^p]$) which exponentially attracts all the
trajectories starting in $v\in D_{+}.$ This is a {\it partial inertial
manifold} for the system (\ref{sdd6-g}), (\ref{sdd6-ic}) with $B_1[\xi]$.


\medskip

As an application we can consider the diffusive Nicholson's
blowflies equation (see e.g. \cite{So-Wu-Yang,So-Yang}) with
state-dependent delay \cite{Rezounenko-Wu-2006,
Rezounenko-JMAA-2007}. More precisely, we consider equation
(\ref{sdd6-g}) where $-A$ is the Laplace operator with the
Dirichlet boundary conditions, $\Omega\subset R^{n_0}$ is a
bounded domain with a smooth boundary,  the nonlinear function $b$
is given by $b(w)=p\cdot w^2e^{-|w|}.$  As a result, we conclude
that under the above assumptions, the diffusive Nicholson's
equation possesses a partial inertial manifold.

\medskip

\noindent {\bf Acknowledgements.}
The author wishes to thank Professor Hans-Otto Walther for bringing
state-dependent delay differential equations to his attention.

\bigskip

\bigskip

June 12, 2007

 Kharkiv

\end{document}